\RequirePackage{pdfmanagement-testphase} 
\DocumentMetadata{testphase} 
\documentclass[12pt]{article}
\tagpdfsetup{tabsorder=structure,uncompress,activate-all,add-new-tag=Title/P,interwordspace=true,tagunmarked=false}  
\usepackage[pdfdisplaydoctitle=true]{hyperref} 
\hypersetup{
    pdftitle={}, 
    pdfauthor={Jeffrey Kuan},
    pdfsubject={},
    pdfkeywords={},
    pdfcreator={Jeffrey Kuan}, 
    bookmarksnumbered=true,     
    bookmarksopen=true, 
    linktocpage=true,
} 
\newcommand{\parag}[1]{#1}
\newcommand{\equ}[2]{#1}

\newcommand{\HeadTwo}[1]{#1}
\newcommand{\HeadThree}[1]{#1}
\newcommand{\HeadFour}[1]{#1}

\usepackage{setspace}
\usepackage[utf8]{inputenc}
\usepackage[T1]{fontenc}
\usepackage{graphicx}
\graphicspath{ {./images/} }
\usepackage{amsmath}
\usepackage{amsfonts}
\usepackage{amssymb}
\usepackage{amsthm}
\theoremstyle{definition}
\newtheorem{defn}{Definition}[section]
\newtheorem{thm}{Theorem}[section]
\newtheorem{prop}{Proposition}[section]

\usepackage{bbold}
\usepackage{fullpage}

\begin{document}

\title{\(q\)--exchangeable Measures and Transformations in Interacting Particle Systems}
\author{Jeffrey Kuan}
\date{}
\maketitle

\begin{center}
\fbox{\parbox{0.8\textwidth}{Accessibility statement: a WCAG 2.1AA compliant version of this PDF will eventually be available at \href{http://tailorswiftbot.godaddysites.com/september-2024-paper}{the author's webpage.}}}\end{center}

\abstract{This paper provides unified calculations regarding certain measures and transformations in interacting particle systems. More specifically, we provide certain general conditions under which an interacting particle system will have a reversible measure, gauge transformation, or ground state transformation. Additionally, we provide a method to prove that these conditions hold. This method uses certain quantum groups, and in that context the general conditions specialize to a \(q\)--exchangeable property.}

\normalsize

\HeadTwo{\section{Introduction}}
\parag{In recent years, there have been developments in constructing interacting particle systems and Markov chains using algebraic machinery. Using certain symmetries arising from algbraic objects, one can construct a stochastic matrix satisfying certain ``nice'' properties, such as Yang--Baxter integrability or Markov duality. The necessary algebraic background is somewhat abstract, and often not presented in a way that's easily digestible to probabilists. In this short set of notes, we provide an exposition of these methods that ``distills'' the necessary probabilistic properties, in a way that is more readable to those without an algebraic background.
}

\parag{More specifically, given a Hamiltonian with an eigenvector, we provide a small set of assumptions which allow the Hamiltonian to be conjugated into a matrix with the sum--to--unity property (this is sometimes called a ``ground state transformation''). Additionally, the conjugation can be explicitly found from the eigenvector. Furthermore, if the Hamiltonian satisfies the Yang--Baxter equation, then a ``gauge transformation'' will result in a matrix also satisfying Yang--Baxter with a sum--to--unity property.  This set of assumptions occur naturally in the context of \(q\)--exchangable measures and Mallows measures. \cite{GO10,GO11,BC24,BB24,Buf20,BN22} This appears to be approximately the minimal set of assumptions necessary; for example, a further of weakening of the assumptions \cite{KZ23} causes the properties to no longer hold.}

\noindent\parag{\textbf{Acknowledgments.} I would like to thank Lafay Augustin for a helpful conversation at the Institute for Pure and Applied Mathematics' workshop on ``Vertex Models: Algebraic and Probabilistic Aspects of Universality.''  I would also like to thank Erik Brodsky and Lillian Stolberg for pointing out a mistake in an earlier version of these notes.}

\HeadTwo{\section{Probabilistic Results}}

\HeadThree{\subsection{Definitions and Assumptions}}

\parag{First we define the state space. Let} \equ{\(\mu=\left(\mu_1, \ldots, \mu_{n}\right)\)}{mu equals mu 1 to mu n} \parag{denote a sequence of \(n\) non--negative integers. We will use}\equ{ \( \vert \mu\vert\)}{absolute value of mu} \parag{to denote} \equ{\( \mu_1+\ldots+\mu_n.\)}{ mu 1 plus dot dot dot plus mu n} \parag{Let \(B_J\) denote the set of all \(mu\) whose absolute value is \(J.\) The notation \(B_J\) is chosen with its algebraic context, where it will be a basis for a vector space. In a mathematical physics or probabilistic setting, each \(\mu\) is sometimes interpreted as a particle configuration at a single lattice sites, with \(\mu_i\) particles of ``type/species/color'' \(i.\)}

\parag{We let the state space be the set}\equ{\[B_{J}\times \cdots \times B_{J}\]}{B J  cross product dot dot dot B J }\parag{so that \(L\) is the number of lattice sites.}\parag{ Let }\equ{\(\mathcal{H}\)}{curly H}\parag{ denote the global Hamiltonian on \(L\) sites, written as:}
\equ{\[\mathcal{H} = H_{12} + H_{23} + \ldots + H_{L-1,L}\]}{curly H equals H 1 2 plus H 2 3 plus dot dot dot H L minus 1 comma L}\parag{where the subscripts \(i,i+1\) indicate that the local Hamiltonian \(H\) acts on lattice sites \(i\) and \(i+1.\)}

 \parag{ Assume that there is an eigenvector \(w\) of \(H\) with eigenvalue \(a,\) or in other words }\equ{\(Hw=aw.\)}{H w equals a w.} \parag{We would like to construct} \equ{\(\mathbf{w}\)}{bold w} \parag{which is an eigenvector of} \equ{\(\mathcal{H}.\)}{ curly H} \parag{To do this, we need to make an assumption on }\equ{\(\mathbf{w}:\)}{bold w}

\HeadFour{\noindent\textbf{Assumption.}}\parag{Assume that \(H\) is weight--preserving in the sense that }\equ{\(\langle \nu,\nu' | H | \eta,\eta'\rangle\)}{bra nu comma nu prime; H; ket eta comma eta prime} is nonzero only if \equ{\(\nu + \nu' = \eta + \eta'.\)}{nu plus nu prime equals eta plus eta prime}  \parag{Suppose that} \equ{\({w}\)}{w} \parag{ has the form}
\equ{\[
{w} = \sum_{\eta,\eta'} W(\eta,\eta')\vert \eta,\eta'\rangle,
\]}{bold w is a sum over eta and eta' of eta comma eta' kets with coefficient upper W of eta comma eta'}\parag{where the function } \equ{\(W\)}{upper W} \parag{ satisfies the property that } \equ{\[ \prod_{k<L+1} W(\eta^k,\eta^{L+1}) \]}{product over k less than L plus 1 of upper W of eta superscript k comma eta superscript L plus 1}\parag{only depends on the value of} \equ{\(\eta^1+\ldots+\eta^L\)}{eta super 1 plus dot dot dot eta super L} \parag{and} \equ{\(\eta^{L+1}.\)}{eta super L plus 1} \parag{Assume also that the quantity}\equ{\[\prod_{k=1}^{L-1} W(\eta^k,\eta^{L})W(\eta^k,\eta^{L+1}) \]}{product from k equals 1 to L minus 1 of: upper W of eta super k comma eta super L; times upper W of eta super k comma eta super L plus 1}\parag{depends only on the values of }\equ{\(\eta^1+\ldots+\eta^{L-1}\)}{eta super 1 plus dot dot dot eta super L minus 1} \parag{and}\equ{ \(\eta^L+\eta^{L+1}\).}{eta super L plus eta super L plus 1}

\HeadFour{\noindent \textbf{Remark.}}\parag{An example of a function \(W\) satisfying this assumption is:}
\equ{
\[
W(\eta,\eta') = q^{2\sum_{i<j} \eta_i \eta_j'}.
\]
}{upper W of eta comma eta prime equals q squared to the power: sum over i less than j of eta sub i times eta sub j prime}

\HeadThree{\subsection{Eigenvectors of Hamiltonians}}

\HeadFour{\begin{thm}}\label{Thm1}
\parag{Make the assumption above and define} \equ{\(\mathbf{w}\)}{bold w} \parag{by}
\equ{\[\mathbf{w} = \sum_{\eta^1,\ldots,\eta^L} \mathbf{W}^L(\eta^1,\ldots,\eta^L)\vert \eta^1,\ldots,\eta^L \rangle,\]}{bold w equals a sum over eta superscript 1 dot dot dot eta superscript L of kets eta superscript 0 dot dot dot eta superscript L. The coefficient is bold upper W superscript L of: eta superscript 0 dot dot dot eta superscript L}\parag{where}\equ{\[\mathbf{W}^L(\eta^1,\ldots,\eta^L)=\prod_{k<l} W(\eta^k,\eta^l)\]}{bold upper W superscript L of eta superscript 1 dot dot dot eta superscript L equals a product over k less than l of upper W of eta superscript k comma eta superscript l }\parag{Then}
\equ{\[ \mathcal{H}\mathbf{w}=(L-1)a\mathbf{w}. \]}{curly H applied to bold w equals the constant L minus 1 times a times bold w}
\end{thm}
\noindent \HeadFour{\textbf{Proof.}} \parag{Proceed by strong induction on \(L.\) The base case \(L=2\) holds by assumption, so now assume the theorem holds for the values from 2 to \(L-1.\)}

\parag{We first re--write the proposed eigenvector} \equ{\(\mathbf{w},\)}{bold w} \parag{as} \equ{\begin{align*}\mathbf{w}&=\sum_{\eta^1,\ldots,\eta^{L+1}}  \mathbf{W}^{L+1}(\eta^1,\ldots,\eta^{L+1})\vert \eta^1, \ldots,\eta^{L+1}\rangle\\ &= \sum_{\eta^1,\ldots,\eta^{L+1}}\left(\prod_{1 \leq i < j \leq L} W(\eta^i,\eta^j)\Big\vert  \eta^1, \ldots,\eta^L \Big\rangle \otimes  \prod_{k<L+1} W(\eta^k,\eta^{L+1})  \Big\vert\eta^{L+1}\Big\rangle\right). \end{align*}}{ }\parag{Since tensor products are bi--linear, we can write this as}
\equ{
\[
\mathbf{w} =\sum_{\eta^1,\ldots,\eta^{L}}\left(\prod_{1 \leq i < j \leq L} W(\eta^i,\eta^j)\Big\vert  \eta^1, \ldots,\eta^L \Big\rangle \otimes \sum_{\eta^{L+1}}  \prod_{k<L+1} W(\eta^k,\eta^{L+1})  \Big\vert\eta^{L+1}\Big\rangle\right)
\]}{}
\parag{Recall that the global Hamiltonian} \equ{\(\mathcal{H}\)}{curly H} \parag{is weight--preserving and that}\equ{\[ \prod_{k<L+1} W(\eta^k,\eta^{L+1}) \]}{product over k less than L plus 1 of upper W of eta superscript k comma eta superscript L plus 1}\parag{only depends on the value of} \equ{\(\eta^1+\ldots+\eta^L\)}{eta super 1 plus dot dot dot eta super L} \parag{and} \equ{\(\eta^{L+1}.\)}{eta super L plus 1} \parag{Since the global Hamiltonian conserves the former quantity, this means that }\equ{\[ (H_{12}+\cdots+H_{L-1,L})\mathbf{w} = (L-1)a\mathbf{w},\]}{H 1 2 plus dot dot dot H L minus 1 L applied to bold w equals the constant L minus 1 times a times bold w}\parag{by the induction hypothesis.}

\parag{So it remains to show that} \equ{\(H_{L,L+1}\mathbf{w} = aw.\)}{H sub L comma L plus 1 applied to bold w equals a times bold w} \parag{This time we write the proposed eigenvector as}\equ{\begin{multline*} \sum_{\eta^L, \eta^{L+1}} \Big( \sum_{\eta^1,\ldots,\eta^{L-1}} \mathbf{W}^{L-1}(\eta^1,\ldots,\eta^{L-1}) \\ \prod_{k=1}^{L-1} W(\eta^k,\eta^{L})W(\eta^k,\eta^{L+1}) \vert \eta^1,\ldots,\eta^{L-1} \rangle \otimes W(\eta^L,\eta^{L+1} )\vert \eta^L,\eta^{L+1} \rangle \Big).\end{multline*}}{}\parag{As before, the quantity}\equ{\[\prod_{k=1}^{L-1} W(\eta^k,\eta^{L})W(\eta^k,\eta^{L+1})\] }{product from k equals 1 to L minus 1 of: upper W of eta super k comma eta super L; times upper W of eta super k comma eta super L plus 1}\parag{depends only on the values of }\equ{\(\eta^1+\ldots+\eta^{L-1}\)}{eta super 1 plus dot dot dot eta super L minus 1} \parag{and}\equ{ \(\eta^L+\eta^{L+1}\).}{eta super L plus eta super L plus 1}\parag{The latter quantity is conserved by the local Hamiltonian, so therefore} \equ{\(H_{L,L+1}\mathbf{w}=a\mathbf{w}.\)}{H sub L comma L plus 1 applied to bold w equals A times bold w.}

\HeadThree{\subsection{Ground State Transformation}}
\parag{ The ground state transformation is actually a special case of Theorem \ref{Thm1}. The definition we use for the ground state transformation is the following:}

\HeadFour{\noindent \textbf{Definition.}} \parag{ Given a Hamiltonian} \equ{\(\mathcal{H}\),}{curly H} \parag{ a \underline{ground state transformation} is a diagonal matrix \(G\) such that all the rows of} \equ{ \(G^{-1}\mathcal{H}G\)}{G inverse curly H G} \parag{ have the same sum.}

\parag{The requirement that all the rows have the same sum is naturally relevant for generators of continuous--time processes and for stochastic matrices.}

\HeadFour{\begin{thm}} \parag{Make the same assumption above, and additionally assume that the function \equ{\(W\)}{upper W} is always nonzero. Then a ground state transformation for the global Hamiltonian } \equ{\(\mathcal{H}\)}{curly H} \parag{always exists, and its entries are given by the entries of} \equ{\(\mathbf{w}.\)}{bold w}
\end{thm}

\HeadFour{\noindent\textbf{Proof.}} \parag{Since} \equ{\(\mathbf{w}\)}{bold w} \parag{ is an eigenvalue of the global Hamiltonian} \equ{\(\mathcal{H},\)}{curly H} \parag{ then for all \(x\) we have}\equ{\[\sum_y \mathcal{H}(x,y)G(y) = \sum_y \mathcal{H}(x,y)\mathbf{w}(y) = (\mathcal{H}\mathbf{w})_x =(L-1)a G(x).\]}{sum over y of curly H of x comma y times G of y; equals same sum with G replaced by bold w; equals the x entry of curly H times bold w; equals L minus 1 times a times G of x}\parag{Dividing by} \equ{\(G(x)\)}{G of x} \parag{completes the proof.}

\HeadThree{\subsection{Gauge Transformation and Yang--Baxter}}

\parag{Given an \(R\)--matrix satisfying the Yang--Baxter equation, we wish to find a ``gauge transformation'' also satisfying the Yang--Baxter equation, such as in \cite{KMMO16}.}

\HeadFour{\begin{thm}}
\parag{Suppose that there is an \(R\)--matrix satisfying the Yang--Baxter equation}\equ{\[ R_{12}R_{13}R_{23} = R_{23}R_{13}R_{12}.\]}{R 1 2 R 1 3 R 2 3 equals R 2 3 R 1 3 R 1 2}\parag{Furthermore, suppose there is a diagonal matrix \(G\) which defines a gauge transformation by}\equ{\[S_{ij}=G_{ji}^{-1}R_{ij}G_{ij}.\]}{ S i j equals G j i inverse times R i j times G i j}\parag{Also assume that for any distinct \(i,j,k,\) there are diagonal matrices \(A_{i,jk},B_{i,jk}\) acting on the \(j,k\) sites such that the following relations hold:}\equ{\begin{align*} G_{ik}^{-1}R_{jk}G_{ik} &= A_{i,jk}R_{jk}A_{i,jk}^{-1}, \\ G_{ij}^{-1}R_{jk}G_{ij} &= A^{-1}_{i,jk}R_{jk}A_{i,jk}, \\  G_{ki}^{-1}R_{jk}G_{ki} &= A^{-1}_{i,jk}R_{jk}A_{i,jk}, \\  G_{ji}^{-1}R_{jk}G_{ji} &= A_{i,jk}R_{jk}A^{-1}_{i,jk} . \end{align*}}{}\parag{Then the \(S\)--matrix also satisfies the Yang--Baxter equation, i.e.}
\equ{\[ S_{12}S_{13}S_{23} = S_{23}S_{13}S_{12}.\]}{S 1 2 S 1 3 S 2 3 equals S 2 3 S 1 3 S 1 2}
\end{thm}

\HeadFour{{\noindent \textbf{Proof.}}}\parag{ We want to prove that}\equ{\[G_{21}^{-1}R_{12}G_{12}G_{31}^{-1}R_{13}G_{13}G_{32}^{-1}R_{23}G_{23} = G_{32}^{-1}R_{23}G_{23}G_{31}^{-1}R_{13}G_{13}G_{21}^{-1}R_{12}G_{12}. \]}{}\parag{Using the relations, we can move the \(G\) terms all the way to the left or all the way to the right. More specifically, on the left--hand--side move \(G_{12}\) and \(G_{13}\) all the way to the right, and move \( G_{31}^{-1}\) and \(G_{32}^{-1}\) all the way to the left. On the right--hand--side, move \(G_{23},G_{13}\) all the way to the right, and move \(G_{31}^{-1},G_{21}^{-1}\) all the way to the left. When making these moves, all the \(A\) terms cancel because:}\equ{\begin{align*} R_{12}G_{31}^{-1}G_{32}^{-1}&=G_{31}^{-1}G_{32}^{-1} R_{12},\\ G_{12}G_{23}R_{23} &= R_{23}G_{12}G_{23} \\ R_{23}G_{31}^{-1}G_{21}^{-1} &= G_{21}^{-1}G_{31}^{-1}R_{23} ,\\ G_{23}G_{13}R_{12} &= R_{12}G_{13}G_{23},\\ G_{12}R_{13}G_{32}^{-1} &= G_{32}^{-1}R_{13}G_{12}\\ G_{23}R_{13}G_{21}^{-1} &=  G_{21}^{-1}R_{13}G_{23}.\end{align*}}{}\parag{After all the moves, the \(G\) terms cancel, and then the result follows because the \(R\)--matrices satisfy the Yang--Baxter equation.} 

\begin{flushright}{\parag{Q.E.D.}}\end{flushright}


\HeadTwo{\section{Algebraic Context}}

\parag{In this section, we provide some algebraic contexts for why the assumptions should hold. Roughly speaking, \(q\)--exchangeable measures occur naturally when there is an underlying algebraic structure. For exclusion processes, \(q\)--exchangeable (or Mallows) measures are defined by the property that swapping adjancent particles multiples the measure by \(q,\) although for partial exclusion processes the definition is more involved \cite{KuanAHP}.}

\parag{In particular, if there are reversible measures which are \(q\)--exchangeable measures, then it follows that the eigenvector \(w\) can be constructed with eigenvalue \(0\). We provide three separate contexts in which these measures can be constructed. We also note that one only needs the measure on the two--lattice process, which avoids more complicated calculations such as the \(q\)--exponential \cite{CGRS,CGRS2,KIMRN}.}

\parag{We introduce some algebraic notation. Let \(V_J\) be a vector space with basis \(B_J.\) Suppose that  }\equ{$\vert \lambda\rangle\in V_J,$}{ket lambda is in V J}\parag{ such that the local Hamiltonian \(H\) satisfies}\equ{\[ H \vert \lambda,\lambda \rangle =0, \qquad S_kH=HS_k, \]}{H applied to the ket lambda comma lambda is zero; and the coproduct of E j commutes with H}\parag{where \(\{S_k\}\) are some operators on } \equ{\(V_J \otimes V_J.\)}{V J tensor V J}\parag{ We will also assume that \(V_J\) is generated by all possible products of \(\{S_k\}\) on}\equ{\(\vert \lambda,\lambda \rangle.\)}{ket lambda comma lambda}\parag{ In other words, \(V_J\) is irreducible over \(\{S_k\}.\) This irreducibility assumption will ensure that the ground state transformation is nonzero, barring miraculous cancelations.}

\HeadThree{\subsection{Quantum Groups}}

\parag{Quantum groups have found use in interacting particle systems, such as \cite{Sch97}, \cite{BS,BS2,BS3}, \cite{CGRS,CGRS2}, \cite{FKZ22,KZ23}, \cite{KLLPZ,REU2022,REU2023}, \cite{KuanJPhysA,KuanCMP,KIMRN,KuanSW}. Here, we briefly explain its relationship to the content of the present paper }

\parag{Here, we define an algebra which is similar to a quantization of }\equ{\(\tilde{D}_{n-1},\)}{tilde D sub n minus 1}\parag{ the type }\equ{\(D\)}{D}\parag{ Lie algebra. We choose this because it is a simply--laced Lie algebra, and contains the other simply--laced Lie algebras as Lie subalgebras. Here we only use the positive Borel subalgebra, to maintain more generality.}

\HeadFour{\begin{defn}}
\parag{Let \(U\) denote the algebra with generators \(E_1,\ldots,E_{n-1},K_1,\ldots,K_{n-1}\) and \(E_1',E_{n-1}',K_1',K_{n-1}'\) and relations }\equ{\[ E_{i}^{\prime} E_{i}^{2}-\left(q+q^{-1}\right) E_{i} E_{i}' E_{i}^{}+E_{i}^{2} E_{i}^{\prime}=0     \]}{E i prime E i squared minus; q plus q inverse times E i E i prime E i; plus E i squared E i prime equals zero}\equ{\[ E_{i\pm 1}^{} E_{i}^{2}-\left(q+q^{-1}\right) E_{i} E_{i\pm 1} E_{i}^{}+E_{i}^{2} E_{i \pm 1}^{}=0     \]}{E i plus of minus 1 E i squared minus; q plus q inverse times E i E i plus or minus 1 E i; plus E i squared E i plus or minus 1 equals zero}\equ{\[K_iE_i = q^2E_iK_i, \quad K_i'E_i'=q^2E_i'K_i', \quad [K_i,K_j]=[K_i,K_j']=[K_i',K_j']=0.\]}{K i E i equals q squared E i K i ; K i prime E i prime equals q squared i prime K i prime; commutator of all the K equals zero.}\equ{\[K_i'E_{i\pm 1}=q^{-1}E_{i \pm 1}K_i', \quad K_iE_{i\pm 1}=q^{-1}E_{i \pm 1}K_i.\]}{K i prime E i plus or minus 1 equals q inverse E i plus or minus 1 K i prime; K i E i plus or minus 1 equals q intersve E i plus or minus 1 K i}
\end{defn}

\HeadFour{\textbf{Remark.}} \parag{ Note that the algebra does not contain the \(F\) elements that occur in a quantum group. These elements would be relevant for establishing the commutation relations \(S_kH=HS_k.\) However, recent work has found other methods for establishing these relations, and therefore do not include the \(F\) elements in the algebra.}%


\parag{We will define a representation, using the ``ket'' notation from mathematical physics. For any non--negative integer \(n,\) define its \(q\)--deformation by}\equ{ \[ [n]_q= \frac{q^n - q^{-n}}{q-q^{-1}}.\]}{fraction with numerator q to the n minus q to the minus n; denominator is q minus q inverse}\parag{As some more notation, let} \equ{\[\epsilon_{i}' =(0, \ldots, 0,1,1,0, \ldots, 0)\]}{epsilon i prime }\parag{consist of a sequence of \(n\) integers, where there is a \(0\) everywhere, except for a \(1\) at the \(i\) and \(i+1\) entries.} \parag{Similarly, let}\equ{\[\epsilon_{i} =(0, \ldots, 0,1,-1,0, \ldots, 0)\]}{epsilon i }\parag{consist of a sequence of \(n\) integers, where there is a \(0\) everywhere, except for a \(1\) at the \(i\) entry and \(-1\) at the \(i+1\) entries.}

\HeadFour{\begin{prop}}
 \parag{The following defines a representation of \(U\):}
\equ{\[
 E_{i}|\mu\rangle=[\mu_{i+1}]_q \left|\mu+\epsilon_{i}\right\rangle \quad   E_{i}^{\prime}|\mu\rangle=|\mu+\epsilon_{i}^{'}\rangle \quad   K_{i}|\mu\rangle=q^{\mu_{i}-\mu_{i+1}}|\mu\rangle \quad K_{i}^{'}|\mu\rangle=q^{\mu_{i}+\mu_{i+1}}|\mu\rangle 
.\]}{E i maps mu to mu plus epsilon i ; times the q deformed mu with the subscript i plus 1; E i prime maps mu to mu plus epsilon i prime; K i multiplies mu by q to the mu sub i minus mu sub i plus 1; K i prime multiplies mu by q to the mu i plus mu i plus 1}
\end{prop}
\HeadFour{\noindent\textbf{Proof.}}\parag{It is a straightforward calculation that this is a representation. For example:}
\equ{\[
\begin{array}{ll}   
  K_{i} E_{i}|\mu\rangle=K_{i}\left|\mu+\epsilon_{i}\right\rangle=q^{\mu_{i}-\mu_{i+1}+2}\left|\mu+\epsilon_{i}\right\rangle=q^{2} E_{i} K_{i}|\mu\rangle \\   K_{i}^{\prime} E_{i}^{\prime}|\mu\rangle=K_{i}^{\prime}\left|\mu+\epsilon_{i}^{\prime}\right\rangle=q^{\mu_{i}+\mu_{i+1}+2}\left|\mu+\epsilon_{i}^{\prime}\right\rangle=q^{2} E_{i}^{\prime} K_{i}^{\prime}|\mu\rangle   \end{array}
\]}{K i E i of mu equals; K i of mu plus epsilon i; equals; q to the power two plus  mu i minus mu sub i plus 1 of mu plus epsilon i; equals q squared E i K i of mu; the second line is the same as the first line but with primes and the minus replaced with a plus; more specifically;  K i prime E i prime of mu equals; K i prime of mu plus epsilon i prime ; equals; q to the power two plus  mu i plus mu sub i plus 1 of mu plus epsilon i; equals q squared E  i prime K i prime of mu; } \parag{Meanwhile, the Serre relations (the first two relations) follow from}
\equ{\[
\left[\mu_{i+1}-1\right]_{q}\left[\mu_{i+1}\right]_{q}-\left(q+q^{-1}\right) \left[\mu_{i+1}\right]_{q}^{2}+\left[\mu_{i+1}\right]_{q}\left[\mu_{i+1}+1\right]_{q} =0  ,
\]}{}\parag{The remaining calculations are similar. Q.E.D. }

\parag{In this setting, the vector}\equ{ \(\vert \lambda\rangle\)}{\lambda}\parag{ is \(0,0,\ldots,J\) and the symmetrices are \(E_1,\ldots,E_{n-1}.\) Some possible ways of showing commutation with the Hamiltonian are with the Casimir \cite{CGRS,CGRS2,KIMRN} or the \(R\)--matrix \cite{KuanCMP}. We also refer to the review paper \cite{CRV14} for relationships between the quantum group and the Zamolodchkov algebra.  }

\HeadThree{\subsection{Hecke algebras}}

\parag{Associated to any Coxeter group is a corresponding Iwahori--Hecke algebra. For simplicity, we consider the type \(A\) Coxeter group, which is simply the symmetric group. Hecke algebras have found use in interacting particle systems \cite{Buf20,KuanSW}.}

\parag{This Hecke algebra is generated by elements \(T_1,\ldots,T_{L-1}\) with relations}\equ{\[ T_iT_{i\pm 1}T_i = T_{i \pm 1}T_i T_{i \pm 1}, \quad T_iT_j=T_jT_i \text{ for } \vert i-j\vert \geq 2, \quad (T_i-q)(T_i+1)=0.\]}{T i times T i plus or minus 1 times T i equals T i plus or minus 1 T i T i plus or minus 1; T i commutes with T j if i and j differ by more than 2; the prodct of T i minus q and T i plus 1 equals 0}\parag{For any dimension \(d,\) there is a representation on}\equ{ \( \underbrace{\mathbb{C}^J \otimes \cdots \otimes \mathbb{C}^J}_{L \text{ tensor products}} \) }{L fold tensor product of d dimensional complex vector space}\parag{ by letting \(T_i\) act on the \(i,i+1\) tensor powers as the \(R\)--matrix}\equ{\[\sum_{i<j} E_{ji}\otimes E_{ij} + q^{-1}\sum_{i<j}  E_{ij}\otimes E_{ji} + (1-q^{1})\sum_{i<j}E_{jj}\otimes E_{ii} + \sum_i E_{ii}\otimes E_{ii}, \]}{}\parag{where \(E_{ij}\) indicates a matrix with a \(1\) at the \(i,j\)--entry and \(0\) elsewhere.}

\parag{Commutation with the Hamiltonian can be shown with Schur--Weyl duality with the quantum groups \cite{KuanSW} or with the Temperley--Lieb algebra, although the latter has not been pursued in a probabilistic setting. }
 
\HeadThree{\subsection{\(q\)--KZ equations}}

\parag{We also mention a relationship of multi--species ASEP to \(q\)--KZ (Knizhnik--Zamolodchikov) relations. This construction produces different algebraic results \cite{CGW20} from the quantum group approach. In this setting, one considers a polynomial--valued vector}\equ{ \(\vert \Phi \rangle \in \mathbf{C}[z_1,\ldots,z_L] \otimes V^{\otimes L}\)}{Phi ket in complex polynomial ring in z 1 to z L tensor product V tensor power L}\parag{ where \(V\) is a vector space. The \(q\)--KZ relations then read}\equ{\[s_i \vert \Phi \rangle = R(z_i/z_{i+1})\vert \Phi \rangle\]}{s sub i applied to Phi equals R of z i over z i plus one applied to Phi}\parag{where \(s_i\) permutes \(z_i\) and \(z_{i+1}\) in the polynomial ring and \(R\) is the \(R\)--matrix. Note that this differs from the presentation of the \(R\)--matrix above, where it is presented as the ``universal'' \(R\)--matrix with no spectral parameters \(z_i.\)}

\parag{It turns out that for \(V=\mathbf{C}^n,\) this can be related to the \( (n-1)\)--species ASEP. More specifically, the \(q\)--KZ relations can be written as} \equ{\(L\vert Phi\rangle = M\vert Phi \rangle\)}{L of phi equals M of phi}\parag{where \(L\) acts on the polynomial ring as a multi--species ASEP generator, and \(M\) acts on \(V^{\otimes L}\) as a multi--species ASEP generator. See \cite{CGW20} for more details. Due to the algebraic background of the \(q\)--KZ relations, it may be possible to extend these results to \(\mathbf{C}^J\) where there is partial exclusion, but this has not yet been pursued in the literature.}

\end{document}